\documentclass[11pt]{article}

\usepackage{amsmath}
\usepackage{float}
\usepackage{amsfonts}
\usepackage{amstext}
\usepackage{graphics}
\usepackage{epsf}
\usepackage[dvips]{epsfig}
\usepackage[dvips]{color}
\usepackage{fullpage}

\def \mod#1 #2{#1\ ({\rm mod}\ #2)}
\def \nodiv{{|\kern-3.9pt/}}

\def \prend{\vrule depth-1pt height7pt width6pt}

\def \endpf{{\ \ \prend \medbreak}}
\newtheorem{theorem}{Theorem}

\newtheorem{conjecture}[theorem]{Conjecture}

\newfont{\cyr}{wncyr10 at 12.0pt}

\def\xx{\phantom{\mbox{\Large $I_I$}}\hspace*{-12pt}}
\def\xxx{\phantom{\mbox{\huge $I_{I_I}$}}\hspace*{-23pt}}
\newcommand{\cnj}[1]{{\textcolor[gray]{0.4}{#1}}} 

\title{A Generalization of Repetition Threshold}

\author{Lucian Ilie\thanks{Research supported in part by NSERC grant R3143A01.}\\
Department of Computer Science\\
University of Western Ontario\\
London, ON,  N6A 5B7 \\
CANADA \\
{\tt ilie@csd.uwo.ca}
\and
Jeffrey Shallit\thanks{Research supported in part by NSERC.}\\
School of Computer Science \\
University of Waterloo \\
Waterloo, ON, N2L 3G1 \\
CANADA \\
{\tt shallit@graceland.uwaterloo.ca}
}

\begin{document}
\date{\today}
\maketitle

\begin{abstract}
     Brandenburg and (implicitly) Dejean
introduced the concept of {\it repetition threshold}:  the smallest
real number $\alpha$ such that there exists an infinite word over a $k$-letter
alphabet that avoids $\beta$-powers for all $\beta > \alpha$.  We
generalize this concept to include the lengths of the avoided words.  We
give some conjectures supported by numerical evidence and prove one of these
conjectures.
\end{abstract}

\section{Introduction}

     In this paper we consider some 
variations on well-known theorems about avoiding repetitions in words.

     A {\it square} is a repetition of the form $xx$, where $x$ is a
nonempty word; an example in English is {\tt hotshots}.  It is easy
to see that every word of length $\geq 4$ over an alphabet of
two letters must contain a square, so squares cannot be avoided in
infinite binary words.  However, Thue showed \cite{Thue:1906,Thue:1912,Berstel:1995}
that there exist infinite words over a three-letter alphabet that avoid
squares.

      Instead of avoiding {\it all\/} squares, one interesting variation
is to avoid 
{\it all sufficiently large\/} squares.  Entringer, Jackson, and Schatz
\cite{Entringer&Jackson&Schatz:1974} showed that there exist infinite
binary words avoiding all squares $xx$ with $|x| \geq 3$.  
Furthermore, they proved that every binary word of length $\geq 18$
contains a factor of the form $xx$ with $|x| \geq 2$, so the bound
$3$ is best possible. 
For some other papers about avoiding sufficiently
large squares, see
\cite{Dekking:1976,Prodinger&Urbanek:1979,Fraenkel&Simpson:1995,Rampersad&Shallit&Wang:2003,Shallit:2003}.

      Another interesting variation is to consider avoiding fractional
powers.  For $\alpha \geq 1 $ a rational number, we say that $y$ is an $\alpha$-power
if we can write $y = x^n x'$ with $x'$ a prefix of $x$ and
$|y| = \alpha |x|$.  For example, the word {\tt alfalfa} is a ${7/3}$-power and
the word {\tt tormentor} is a ${3 \over 2}$-power.
For real $\alpha > 1$, we
say a word {\it avoids $\alpha$-powers} if it contains no factor that is
a $\alpha'$-power for any rational $\alpha' \geq \alpha$.  Brandenburg \cite{Brandenburg:1983} and
(implicitly) Dejean \cite{Dejean:1972} considered
the problem of determining the {\it repetition threshold};
that is, the least exponent $\alpha = \alpha(k)$ such that there
exist infinite words over an alphabet of size $k$ that
avoid $(\alpha+ \epsilon)$-powers for all $\epsilon > 0$.  Dejean
proved that $\alpha(3) = 7/4$.  She also conjectured that $\alpha(4) = 7/5$ and
$\alpha(k) = k/(k-1)$ for $k \geq 5$.  In its full generality, this conjecture is still open,
although Pansiot \cite{Pansiot:1984c} proved that $\alpha(4) = 7/5$ and
Moulin-Ollagnier \cite{Moulin-Ollagnier:1992} proved that Dejean's
conjecture holds for $5 \leq k \leq 11$.  For more information,
see \cite{Choffrut&Karhumaki:1997}.

     In this paper we consider combining these two variations.  We generalize
the repetition threshold of Dejean to handle
avoidance of all sufficiently large fractional powers.
(Pansiot also suggested looking at this generalization
at the end of his paper \cite{Pansiot:1984c}, but to the best of our
knowledge no one else has pursued this question.) We 
give a large number of conjectures, supported by numerical evidence,
about generalized repetition threshold, and prove one of them.

\section{Definitions}

      Let $\alpha > 1$ be a rational number,
and let $\ell \geq 1$ be an integer.
A word $w$ is a {\it repetition of order $\alpha$ and length $\ell$} if
we can write it as $w = x^n x'$ where $x'$ is a prefix of $x$,
$|x| = l$, and $|w| = \alpha |x|$.  For brevity, we also call $w$ a
{\it $(\alpha, \ell)$-repetition}.
Notice that an $\alpha$-power
is an $(\alpha, \ell)$-repetition for some $\ell$.
We say a word is {\it $(\alpha,\ell)$-free}
if it contains no factor that is a $(\alpha',\ell')$-repetition for $\alpha' \geq \alpha$
and $\ell' \geq \ell$.
We say a word is {\it $(\alpha^+,\ell)$-free}
if it is $(\alpha',\ell)$-free for all $\alpha' > \alpha$.

      For integers $k \geq 2$ and $\ell \geq 1$,
we define the {\it generalized repetition threshold} $R(k,\ell)$
as the real number $\alpha$ such that either

\begin{itemize}
\item[(a)] over a $k$-letter alphabet there exists an $(\alpha^+,\ell)$-free
infinite word, but all $(\alpha,\ell)$-free words are finite; or

\item[(b)] over a $k$-letter alphabet there exists a $(\alpha,\ell)$-free
infinite word, but for all $\epsilon > 0$,
all $(\alpha-\epsilon,\ell)$-free words are finite.
\end{itemize}

      Notice that $R(k,1)$ is essentially the repetition threshold of Dejean
and Brandenburg.

\begin{theorem}
The generalized repetition threshold $R(k,\ell)$ exists and is finite
for all integers $k \geq 2$ and $\ell \geq 1$.  Furthermore,
$1 + \ell/k^\ell \leq R(k,\ell) \leq 2$.
\end{theorem}

\begin{proof}
      Define $S$ to be the set of all real numbers $\alpha\geq 1$ such that
there exists a $(\alpha,\ell)$-free infinite word
over a $k$-letter alphabet.   Since Thue proved that there exists
an infinite word over a two-letter alphabet (and hence over larger
alphabets) avoiding $(2+\epsilon)$-powers for all $\epsilon > 0$,
we have that $\beta = \inf S$ exists and $\beta \leq 2$.
If $\beta \in S$, we are in case (b) above, and if $\beta \not\in S$, we
are in case (a).  Thus $R(k,\ell) = \beta$.  

      For the lower bound, note that any word of length $\geq k^\ell + \ell$
contains $\geq k^\ell + 1$ factors of length $\ell$.  Since there are only
$k^\ell$ distinct factors of length $\ell$, such a word
contains at least two occurrences of some word of length $\ell$, and hence
is not $(1+ {\ell\over{k^\ell}}, \ell)$-free.
\endpf
\end{proof}

\bigskip

\noindent{\bf Remarks.}

1. It may be worth noting that we know no instance where case (b) of the definition of
generalized repetition threshold above actually
occurs, but we have not been able to rule it out.

2.  Using the Lov\'asz local lemma,
Beck \cite{Beck:1981} has proved a related result:  namely, for
all $\epsilon > 0$, there exists an integer $n'$ and an infinite
$(1 + n/(2-\epsilon)^n, n)$-free binary word for all $n \geq n'$.
Thus our work can be viewed as a first attempt at an explicit version of Beck's result
(although in our case the exponent does not vary with $n$).

\section{Conjectures}

     In this section we give some conjectures about $R(k,\ell)$.

\begin{conjecture}
\begin{displaymath}
     R(k,2) = \begin{cases}
	(3k-2)/(3k-4), & \text{if $k$ is even;} \\
	(3k-3)/(3k-5), & \text{if $k$ is odd.}
\end{cases}
\end{displaymath}
\end{conjecture}

\begin{conjecture}
     $R(3,\ell) = 1 + {1\over\ell}$ for $\ell \geq 2$.
\end{conjecture}

      These conjectures are weakly supported by some numerical evidence.
The following table
gives the established and conjectured values of $R(k,\ell)$.
Entries in {\bf bold} have
been proved; other entries, in light gray, are merely conjectured.
If the entry for $(k, \ell)$ is $\alpha$, then we have proved,
using the usual tree-traversal technique discussed below,
that there is no infinite $(\alpha, \ell)$-free
word over a $k$-letter alphabet.  

\begin{figure}[H]
\begin{center}
$\begin{array}{c|c}
R(k,\ell) & \ell\\
\hline
k & 
\begin{array}{c||ccccccccccc}
\ \ &\ \ \xx 1\ \ &\ \ 2\ \ &\ \ 3\ \ &\ \ 4\ \ &\ \ 5\ \ &\ \ 6\ \ &\ \ 7\ \  \\
% \ \ 8\ \ &\ \  9\ \ &\ \ 10 \ \ & \ \ 11 \ \ \\
\hline\hline
\xxx 2 & \mathbf{2} &\mathbf{2}&\cnj{\frac{8}{5}}&\cnj{\frac{3}{2}}&\cnj{\frac{7}{5}}&\cnj{\frac{4}{3}}&\cnj{\frac{9}{7}} \\
% \cnj{\frac{10}{8}}&\cnj{\frac{11}{9}} &\cnj{\frac{12}{10}}&\cnj{\frac{13}{11}} \\
\xxx 3 & \mathbf{\frac{7}{4}} &\mathbf{\frac{3}{2}}&\cnj{\frac{4}{3}}&\cnj{\frac{5}{4}}&\cnj{\frac{6}{5}}&\cnj{\frac{7}{6}}&\cnj{\frac{8}{7}} \\
%	&\cnj{\frac{9}{8}}\\
\xxx 4 & \mathbf{\frac{7}{5}} &\cnj{\frac{5}{4}}&\cnj{\frac{6}{5}}&\cnj{\frac{7}{6}}\\
% &\cnj{\frac{8}{7}} &\cnj{\frac{9}{8}}&&\\
\xxx 5 & \mathbf{\frac{5}{4}} &\cnj{\frac{6}{5}}&\cnj{\frac{8}{7}}&&&&&\\
\xxx 6 & \mathbf{\frac{6}{5}} &\cnj{\frac{8}{7}}&&&&&&\\
\xxx 7 & \mathbf{\frac{7}{6}} &&&&&&&\\
\xxx 8 & \mathbf{\frac{8}{7}} &&&&&&&\\
\xxx 9 & \mathbf{\frac{9}{8}} &&&&&&&\\
\xxx 10 & \mathbf{\frac{10}{9}} &&&&&&&\\
\xxx 11 & \mathbf{\frac{11}{10}} &&&&&&&\\
\xxx 12 & \cnj{\frac{12}{11}} &&&&&&&\\
\xxx 13 & \cnj{\frac{13}{12}} &&&&&&&\\
\end{array}
\end{array}$
\caption{Known and conjectured values of $R(k,\ell)$.}
\label{rvalues}
\end{center}
\end{figure}

     The proved results are as follows:

\begin{itemize}
\item $R(2,2) = 1$ follows from Thue's proof of the existence of
overlap-free words over a two-letter alphabet \cite{Thue:1906,Thue:1912,Berstel:1995};

\item $R(2,2) = 2$ follows from Thue's proof together with the observation of
Entringer, Jackson and Schatz
\cite{Entringer&Jackson&Schatz:1974};

\item $R(3,1) = 7/4$ is due to Dejean \cite{Dejean:1972};

\item $R(4,1) = 7/5$ is due to Pansiot \cite{Pansiot:1984c};

\item $R(k,1) = k/(k-1)$ for $5 \leq k \leq 11$ is due to Moulin-Ollagnier
\cite{Moulin-Ollagnier:1992};

\item $R(3,2) = 3/2$ is new and is proved in Section~\ref{newsec}.
\end{itemize}

     We now explain how the conjectured results were obtained.  We used the
usual tree-traversal technique, as follows:  suppose we want to determine
if there are only finitely many words over the
alphabet $\Sigma$ that avoid a certain set of
words $S$.  We construct a certain tree $T$ and traverse it using
breadth-first search.  The tree $T$ is defined as follows:  the root
is labeled $\epsilon$ (the empty word).  If a node $w$ has a factor
contained in $S$, then it is a leaf.  Otherwise, it has children
labeled $wa$ for all $a \in \Sigma$.  It is easy to see that $T$
is finite if and only if there are finitely many words
avoiding $S$.  

      We can take advantage of various symmetries in $S$.  For example,
if $S$ is closed under renaming of the letters (as is the case in the
examples we study), we can label the root with an arbitrary single letter 
(instead
of $\epsilon$) and deduce the number of leaves in the full tree by multiplying
by $|\Sigma|$.

      If the tree is finite, then certain parameters about the tree give useful
information about the set of finite words avoiding $S$:

\begin{itemize}

\item If $h$ is the height of the tree, then any word of length 
$\geq h$ over $\Sigma$ contains a factor in $S$.

\item If $M$ is the length of longest word avoiding $S$, then
$M = h-1$.

\item If $I$ is the number of internal nodes, then
there are exactly $I$ finite words avoiding $S$.  Furthermore,
if $L$ is the number of leaves, then (as usual),
$L = 1 + (|\Sigma|-1) I$.

\item If $I'$ is the number of internal nodes at depth $h-1$,
then there are $I'$ words of maximum length avoiding $S$.

\end{itemize}

Table 2 gives the value of some of these parameters.
Here $\alpha$ is the established or conjectured value of
$R(k,\ell)$ from Table 1.

     We have seen how to prove computationally that only finitely many
$(\alpha,\ell)$-free words exist.  But what is the evidence that suggests
we have determined the smallest possible $\alpha$?  For this, we
explore the tree corresponding to avoiding $(\alpha^+, \ell)$-repetitions
using {\it depth-first} (and not breadth-first) search.
If we are able to construct a ``very long'' word avoiding
$(\alpha^+, \ell)$-repetitions, then we suspect we have found the optimal
value of $\alpha$.  For each unproven $\alpha$ given in Table 1, we were
able to construct a word of length at least $500$ (and in some cases, $1000$)
avoiding the corresponding
repetitions.  This constitutes weak evidence of the correctness of our
conjectures, but it is evidently not conclusive.  

     To show what can go wrong, the data we presented evidently suggests
the conjecture $R(2,\ell) = (\ell+2)/\ell$.
But we have proven this is not true, since the tree avoiding
%$({10 \over 8}^+,8)$-repetitions is finite 
%(height 177, with 25295779 internal nodes).
$(1.2608,8)$-repetitions is finite, with height $195$ and
$53699993$ internal nodes.  (Perhaps $R(2,8) = 29/23$.)

\begin{figure}[H]
\begin{center}
\begin{tabular}{|r|r|c|r|r|r|r|r|l|}
\hline
$k$ & $\ell$ & $\alpha$ & $L$ & $I$ & $h$ & $M =$ & $I'$ & lexicographically least word of \\
    &        &          &     &     &     & $h-1$ &      & length $M$ avoiding $(\alpha,l)$-repetitions \\
\hline
2 & 1 & $\mathbf{2}$ & 8 & 7 & 4 & 3 & 2 & {\tiny 010} \\
2 & 2 & $\mathbf{2}$ & 478 & 477 & 19 & 18 & 2 & {\tiny 010011000111001101}  \\
2 & 3 & $8/5$ & 5196 & 5195 & 34 & 33 & 12 & 
	{\tiny 001100001010111100001110101000110}  \\
2 & 4 & $3/2$ & 13680 & 13679 & 54 & 53 & 4 & 
	{\tiny 01110010010111100000110110100100111110000010110110001} \\
2 & 5 & $7/5$ & 40642 & 40641 & 60 & 59 & 4 & 
	{\tiny 00111010101000001111110010001011101100000011111010101000001} \\
2 & 6 & $4/3$ & 21476 & 21475 & 40 & 39 & 4 & 
	{\tiny 000110101101000000011111110101001000110} \\
2 & 7 & $9/7$ & 81368 & 81367 & 65 & 64 & 4 & 
	{\tiny 0001111011100000001010101011111111001001001011011011} \\
  &   &       &       &       &    &    &   &
	\ \ {\tiny 000000001010}
	\\
%2 & 8 & $5/4$ & 73408 & 73407 & 49 & 48 & 2 & 
%	{\tiny 010110111011000000000101010111111111001000100101}
%	\\
%2 & 9 & $11/9$ & 170720 & 170719 & 57 & 56 & 8 & 
%	{\tiny 01010010010001111111111010101010100000000001101110111010} \\
%2 & 10 & $6/5$ & 252946 & 252945 & 42 & 41 & 16 & 
%	{\tiny 01011111111111000000000001010110111011110} \\
%2 & 11 & $13/11$ & 532428 & 532427 & 45 & 44 & 26 & 
%	{\tiny 01010010000000001111111111110101010101010001} \\
\hline\hline
3 & 1 & $\mathbf{7/4}$ & 6393 & 3196 & 39 & 38 & 18 & 
	{\tiny 01020121021201021012021020121021201020}\\
3 & 2 & $\mathbf{3/2}$ & 11655 & 5827 & 31 & 30 & 6 &
	{\tiny 012002112201100221120011022012}  \\
3 & 3 & $4/3$ & 4037361 & 2018680 & 228 & 227 & 6 &
 {\tiny 012121000111222010121200022210102021112220001212020111000} \\
  &   &       &         &         &     &     &   & 
 \ \ {\tiny 212101022200011120201012221110202121000111222010121200022 } \\
  &   &       &         &         &     &     &   & 
 \ \ {\tiny 211120201012220001110202121000222010121200011122210102021} \\
  &   &       &         &         &     &     &   & 
 \ \ {\tiny 11000121202011122200021210102221112020101222000111020201} \\
3 & 4 & $5/4$ & 188247 & 94123 & 63 & 62 & 24 &
	{\tiny 001022021111000012212102002011112222100101202202111} \\
  &   &       &        &       &    &    &    &
	\ \ {\tiny 00001212210} \\
3 & 5 & $6/5$ & 493653 & 246826 & 63 & 62 & 12 &  
	{\tiny 01011121200000222221110102020212121000001111122022002} \\
  &   &       &        &       &    &    &    &
	\ \ {\tiny 101210120} \\
3 & 6 & $7/6$ & 782931 & 391465 & 60 & 59 & 24 &
	 {\tiny 000012112121020202201111110000002122121201010110}\\
  &   &       &        &       &    &    &    &
	 \ \ {\tiny 22222200001} \\
3 & 7 & $8/7$ &  2881125 & 1440562 & 68 & 67 & 24 &
	{\tiny 00001111111222020201010101212121200000002222222110} \\
  &   &       &        &       &    &    &    &
	\ \ {\tiny 11010012020212021} \\
%3 & 8 & $9/8$ &  6987903 & 3493951 & 62 & 61 & 24 &
%	{\tiny 00011111111022020202012121211210000000022222222011} \\
%  &   &       &        &       &    &    &    &
%	\ \ {\tiny 01010102121 } \\
\hline\hline
4 & 1 & $\mathbf{7/5}$ & 709036 & 236345 & 122 & 121 & 48 &
{\tiny 012031021301231032013021031230132031021301203210231201302} \\
  &   &       &        &        &     &     &    &
\ \ {\tiny 1032012310213203123013210231203213012310320130210312301} \\
  &   &       &        &        &     &     &    &
\ \ {\tiny 320310230 } \\
4 & 2 & $5/4$ & 10324 & 3441 & 17 & 16 & 24 & {\tiny 0112330022110332} \\
4 & 3 & $6/5$ & 153724 & 51241 & 24 & 23 & 96 & {\tiny 01012333000222111332001} \\
4 & 4 & $7/6$ & 2501620 & 833873 & 35 & 34 & 24 & {\tiny 0010122223033111100002212333301011} \\
4 & 5 & $8/7$ & 30669148 & 10223049 & 40 & 39 & 864 &
	{\tiny  001012222230331111100000221233333010101 } \\
\hline\hline
5 & 1 & $\mathbf{5/4}$ & 1785 & 446 & 7 & 6 & 120 & {\tiny 012340} \\
5 & 2 & $6/5$ & 453965 & 113491 & 23 & 22 & 240 &
	{\tiny 0122344002114332204413} \\
5 & 3 & $8/7$ & 7497345 & 1874336 & 34 & 33 & 720 &
	{\tiny 010123234440002111433322204041312} \\
\hline\hline
6 & 1 & $\mathbf{6/5}$ & 13386 & 2677 & 8 & 7 & 720 & {\tiny 0123450} \\
6 & 2 & $8/7$ & 3159066 & 631813 & 21 & 20 & 1440 & {\tiny 01233455002211443052} \\
\hline\hline
7 & 1 & $\mathbf{7/6}$ & 112441 & 18740 & 9 & 8 & 5040 & {\tiny 01234560} \\
\hline\hline
8 & 1 & $\mathbf{8/7}$ & 1049448 & 149921 & 10 & 9 & 40320 & {\tiny 012345670} \\
\hline
\end{tabular}
\end{center}
\caption{Tree statistics for various values of $k$ and $l$}
\end{figure}

\section{A new result}
\label{newsec}

      In this section we prove the following new result:

\begin{theorem}
      $R(3,2) = {3 \over 2}$.
\end{theorem}

\begin{proof}
From the numerical
results reported in Table 2, we know that there exist no infinite words
over a $3$-letter alphabet avoiding $({3 \over 2},2)$-repetitions.  It therefore
suffices to exhibit an infinite word over a $3$-letter alphabet that
avoids $({3\over 2}^+,2)$-repetitions.

   Now consider the uniform morphism 
$h : \{0,1,2,3\}^{\star} \longrightarrow \{0,1,2\}^{\star}$ defined by
$$
\begin{array}{lll}
h(0) = 000211, & \ \ \ & h(2) = 020011,\\
h(1) = 101221, & \ \ \ & h(3) = 120221.\\
\end{array}
$$
By a result of \cite{Pansiot:1984c},
there exist $\frac{7}{5}^+$-free infinite words over four letters.
Consider one such word
${\bf x}$.   We will prove that $h({\bf x})$ is $(\frac{3}{2},2)^+$-free. 

We notice first the following synchronising property of the morphism $h$:
for any $a,b,c \in \{0,1,2,3\}$ 
and $s,r\in\{0,1,2\}^{\star}$, if $h(ab) = rh(c)s$,
then either $r=\varepsilon$ and $a=c$ or $s=\varepsilon$ and $b=c$. 
This is straightforward to verify.

We now argue by contradiction. Assume $h({\bf x})$ has a repetition $xyx$ such that $|x| > |y|$.
If $|x|\ge 11$, then each occurrence of $x$ contains as factor at least one full $h$-image of a letter.
By the above synchronising property, the second $x$ will contain the same full images and at the same positions,
say $x = x'x''x'''$ with $x'' = h(u)$, $|x'|\le 5, |x'''|\le 5$.
Therefore, $h({\bf x})$ contains the factor 
$x'h(u)x'''yx'h(u)x'''$ and ${\bf x}$ has the factor $uvu$,
where $h(v)=x'''yx'$. 
We next compute the order of this repetition in ${\bf x}$. 
Assuming $|x|\ge 50$, we have $|x| \ge 5|x'x'''|$ and so
$$
\frac{|uvu|}{|uv|} 
= 1 + \frac{|x''|}{|x|+|y|} =
1 + \frac{|x|-|x'x'''|}{|x|+|y|} > \frac{7}{5},
$$
a contradiction, since ${\bf x}$ is $\frac{7}{5}^+$-free.
For the case when $|x|<50$, 
it can be shown by exhaustive search that the only possibility for such a 
repetition in $h({\bf x})$
is $|xy|=1$. Thus, $h({\bf x})$ is $(\frac{3}{2},2)^+$-free.
\endpf
\end{proof}

\noindent{\bf Remarks.}

1.  We needed much less than $|x|>|y|$ in obtaining 
the contradiction. In fact, $\frac{|x|}{|y|} > \frac{2}{3}+\delta$,
for some $\delta>0$ is sufficient.
What we obtain is, for any $\delta>0$ there exists $k=k(\delta)$ 
such that $h({\bf x})$ is
$(\frac{7}{5}+\delta,k(\delta))$-free. 

2.  The set of those $x,y$ for which we need to check 
$(\frac{3}{2},2)^+$-freeness 
can be substantially reduced by a more detailed analysis. (For instance,
we bounded $|x'x'''|$ by 10, the simplest bound,
but this can be significantly reduced.) We used
$|x|<50$ in order to simplify the proof.

\section{Acknowledgments}

    We thank Narad Rampersad for a careful reading of this paper and
for pointing out the remark of Pansiot.

\newcommand{\noopsort}[1]{} \newcommand{\singleletter}[1]{#1}

\end{document}